\documentclass[11pt,bibtex]{article}
\newtheorem{theorem}{Theorem}
\begin{document}
\title{Pattern avoidance in compositions and multiset permutations}
\author{Carla D. Savage
\thanks{Research supported in part by NSF grant DMS-0300034}\\
North Carolina State University\\
Raleigh, NC 27695-8206\and
Herbert S. Wilf\\
University of Pennsylvania\\
Philadelphia, PA 19104-6395}
\maketitle

\section{Introduction}
One of the most arresting phenomena in the theory of pattern
avoidance by permutations is the fact that the number of
permutations of $n$ letters that avoid a pattern $\pi$ of 3 letters
is independent of $\pi$. In this note we exhibit two generalizations
of this fact, to ordered partitions, a.k.a. compositions, of an
integer, and to permutations of multisets. It is remarkable that the
conclusions are in those cases identical to those of the original
case.

Further, the number of permutations of a multiset $S=1^{a_1}2^{a_2}\dots k^{a_k}$ that avoid a given pattern $\pi\in S_3$ is a symmetric function of the $a_i$'s, and we will give here a bijective proof of this fact for $\pi = (123)$.

By a composition of an integer $n$ into $k$ parts we mean an integer
representation
\[n=x_1+x_2+\dots+x_k\qquad (\forall i:x_i\ge 0)\]
where two compositions are regarded as distinct even if they differ
only in their order of the summands. If in fact we have $x_i\ge 1$
for all $i$ then we speak of a composition into positive parts. A
composition is said to contain the pattern, e.g., $\pi=(132)$ if
$\exists \,i_1<i_2<i_3$ such that $x_{i_1}<x_{i_3}<x_{i_2}$, and
similarly for other patterns $\pi$. Note the strict inequalities
that we use, which, in view of the repeated part sizes that can
occur, are quite material, though of course variants of these
problems can also be considered in which some of the inequalities
might not be strict.

\begin{theorem}
\label{th:th1} Among the $2^{n-1}$ compositions of $n$ into positive
parts, the number that avoid a given pattern $\pi$ of 3 letters is
independent of $\pi$.
\end{theorem}
These numbers, which play a role analogous to that of the Catalan
numbers in the case of permutations, are, for $n=1,2,\dots$,
\[1,2,4,8,16,31,60,114,214,398,732,1334,2410,\dots .\]
 Below we will find the ordinary power series generating
function of the above sequence, in the form
\begin{equation}
\label{eq:gf1} f(x)=\sum_{i\ge 1}\frac{1}{1-x^i}\prod_{j\neq
i}\left\{\frac{1-x^i}{(1-x^{j-i})(1-x^i-x^j)}\right\}.
\end{equation}

A more refined version of Theorem \ref{th:th1}  is also true.
\begin{theorem}
\label{th:th2} Among the ${n-1\choose k-1}$ compositions of $n$ into
$k$ positive parts, the number that avoid a given pattern $\pi$ of 3
letters is independent of $\pi$. Likewise, among the ${n+k-1\choose
n}$ compositions of $n$ into $k$ nonnegative parts the same
conclusion holds.
\end{theorem}

Finally, at the root of all of the above is the following finer
gradation.
\begin{theorem}
\label{th:th3} Fix a multiset $S$. The number of permutations of $S$ that avoid $\pi$ is independent of the choice of $\pi\in S_3$.
\end{theorem}

The above results are all easy consequences of the enumeration of
the permutations of a \textit{multiset} that avoid a pattern. This
enumeration was accomplished for the pattern $(132)$ in \cite{aaa, alw} by an elegant
recursive construction of the sequence of generating functions
involved. Those authors found that the number of
permutations of the multiset $S=1^{a_1}2^{a_2}\dots k^{a_k}$ that avoid the pattern $(132)$ is the coefficient of $x_1^{a_1}\dots x_k^{a_k}$ in the generating function
\begin{equation}
\label{eq:maingf}
g_k(\mathbf{x})=\sum_{i=1}^k\frac{x_i^{k-1}(1-x_i)^{k-2}}{\prod_{1\le
j\le k;\,j\neq i}\left\{(x_i-x_j)(1-x_i-x_j)\right\}}.
\end{equation}
We remark that, in $g_k$, for fixed $i<j$, the coefficient of $1/(x_i-x_j)$ is a skew-symmetric function of $x_i,x_j$, which therefore contains a
factor of $x_i-x_j$ to cancel those factors that seem to appear in
the denominators. Hence $g_k$ is an analytic function of
$\mathbf{x}$ in a neighborhood of the origin in $R^k$. For example,
\[g_2(\mathbf{x})=\frac{1}{1-x_1-x_2},\]
and
\[g_3(\mathbf{x})=\frac{1-x_1-x_2-x_3+x_1x_2+x_1x_3+x_2x_3}{(1-x_1-x_2)(1-x_1-x_3)(1-x_2-x_3)}.\]

\section{Proofs}
If we assume, for a moment, that Theorem \ref{th:th3} has been
proved then Theorems \ref{th:th1} and \ref{th:th2} follow easily.
Indeed, in each case define an equivalence relation on the set of
compositions involved by declaring that two compositions are
equivalent if they have the same multisets of parts. Then Theorem
\ref{th:th3} implies that the desired pattern-independence
conclusion holds separately within each equivalence class, so
\textit{a fortiori} it holds on the full set of compositions being
considered.

The proof of Theorem \ref{th:th3} is an easy extension of the results in \cite{aaa} and \cite{alw}. Let
$\mathbf{a}=(a_1,\dots,a_k)$ be a given vector of $k$ positive
integers. The multiset $M(\mathbf{a})$ is the one that contains
exactly $a_i$ copies of the letter $i$, for each $i=1,\dots,k$.
Define, for each pattern $\pi$, $f(\mathbf{a},\pi)$ to be the number
of permutations of the multiset $M(\mathbf{a})$ that avoid the
pattern $\pi$. Then, in \cite{aaa} and \cite{alw} it was shown
that
\begin{enumerate}
\item \label{it:i}For every $\mathbf{a}$ we have
$f(\mathbf{a},(123))=f(\mathbf{a},(132))$, and
\item \label{it:ii}$f(\mathbf{a},(132))$ is a symmetric function of
$(a_1,\dots,a_k)$. This follows easily from (\ref{eq:maingf}) above. However it seems to be quite a remarkable fact, in that one might imagine that the number of permutations that avoid a pattern might change drastically if we swap the available supplies of large and small letters. We give a bijective proof of this
 symmetry  for the pattern (123) in Section \ref{sec:symm} below.
\end{enumerate}

Now the reversal map shows that $f(\mathbf{a},\pi)$ is constant on
each of the pairs
\[((123),(321)), ((132),(231)), ((213),(312))\]
 of patterns of length 3. Further, item \ref{it:i} above shows that the
first two of these three classes coincide, so it suffices to show
that for all $\mathbf{a}$ we have
$f(\mathbf{a},(132))=f(\mathbf{a},(312))$.

The complement $\sigma'$ of a permutation $\sigma$ of $n$
letters is obtained by replacing each letter $i$ in the values of
$\sigma$ by $n+1-i$. Evidently we have
$f(\mathbf{a},\sigma)=f(\mathbf{\bar{a}},\sigma')$, where
$\bar{a}$ is the reversal of the vector $\mathbf{a}$. Thus
\[f(\mathbf{a},(312))=f(\mathbf{\bar{a}},(132)).\]
Since all four patterns $(123),(321),(132),(231)$ have the same
counting function, it follows from item \ref{it:ii} above that all
four counting functions are symmetric functions of the $a_i$'s. Thus
$f(\mathbf{\bar{a}},(132))=f(\mathbf{a},(132))$, completing the
proof of Theorem \ref{th:th3}.

A proof of Theorem \ref{th:th3} that is independent of the results
of \cite{aaa} and \cite{alw}, relying instead on a generalization of
the methods of Simion and Schmidt \cite{ss}, has been given by Amy
Myers \cite{anm}.

\section{Counting compositions that avoid a pattern}
Let $\pi$ be a fixed pattern of three letters. We will specialize
the generating function (\ref{eq:maingf}) of the previous section
to the enumeration
of compositions of the integer $n$, into positive parts, that avoid
the pattern $\pi$. By Theorem \ref{th:th1}, of course, these results
will be valid for all $\pi\in S_3$.

In the generating function (\ref{eq:maingf}), put $x_i=x^i$, for
$i=1,\dots,k$, and consider the coefficient of $x^n$ in
$g_k(x,x^2,\dots,x^k)$. This will be the sum of the numbers of
multiset permutations that avoid $\pi$, summed over all multisets
for which $a_1+2a_2+\dots+ka_k=n$, i.e., over all compositions of
$n$ into positive parts. This completes the proof of eq.
(\ref{eq:gf1}).

Doron Zeilberger has remarked that since the generating function has infinitely many singularities, the sequence that it generates is not $P$-recursive. Hence pattern avoidance in compositions of an integer provides a simple and natural example of an avoidance problem whose solution sequence is not $P$-recursive.

If, in the generating function (\ref{eq:gf1}) we halt the outer sum
at $i=k$ then we will be looking at the generating function for the
$\pi$-avoiding compositions of $n$ whose parts are all $\le k$.

Now we can deal with the asymptotic growth rate of $c(n,k)$, the
number of compositions of $n$, into positive parts $\le k$, that
avoid the pattern $\pi$. Since the generating function
(\ref{eq:gf1}), with the upper limit of the outer $i$-sum replaced
by $k$, is a rational function, what we need to do is to determine
which one of the trinomial equations
\[x^i+x^j=1\qquad (1\le i<j\le k)\]
has a root closest to the origin in the complex plane.

We claim that the ``golden ratio'' equation $x+x^2=1$ is the winner
of this competition. Indeed, suppose $x$ is a root of $x^i+x^j=1$,
in which $i<j\le k$ and $j>2$. Since the product of the roots is
$-1$, we can suppose $|x|<1$. Then we have
\[1=|x^i+x^j|\le |x|^i+|x|^j<|x|+|x|^2,\]
whence $|x|>2/(1+\sqrt{5})$, as claimed. It follows that if $c(n,k)$
is the number of $\pi$-avoiding compositions of $n$ into positive
parts $\le k$, then
\begin{equation}
\label{eq:casm}
 c(n,k)\sim
K(k)\left(\frac{1+\sqrt{5}}{2}\right)^n\qquad(n\to\infty).
\end{equation}
where
\[K(k)=\frac{r}{(r-1)(r-s)}\left(r\prod_{j = 3}^{k}
        \frac{1 - \frac{1}{r}}
         {\left( 1 - r^{1 - j} \right) \,\left( 1 - \frac{1}{r} - r^{-j} \right) }- \prod_{j = 3}^{k}\frac{1 - r^{-2}}
         {\left( 1 - r^{2 - j} \right) \,\left( 1 - r^{-2} - r^{-j} \right) }           \right),
\]
and $r=(1+\sqrt{5})/2$, $s=(1-\sqrt{5})/2$. Numerically, the values of $K(k)$ for $k=5,10,20,$ and $\infty$ are 9.95025, 17.9099, 18.9314, and 18.9399867..

\section{A bijection to show symmetry}
\label{sec:symm}

As in Section 2,
for
a given vector
$\mathbf{a}=(a_1,\dots,a_k)$
 of $k$ positive integers,
let $M(\mathbf{a})$ be the multiset  containing
exactly $a_i$ copies of the letter $i$, for each $i=1,\dots,k$,
and let
$f(\mathbf{a},\pi)$ be the number
of permutations of $M(\mathbf{a})$ that avoid the
pattern $\pi$.

In this section we give an explicit bijection to
show that  $f(\mathbf{a},(123))$ is a symmetric function of
$(a_1, \ldots, a_k)$.  For any permutation $\mathbf{b}$ of $\mathbf{a}$,
we show there is a bijection
\[
\Theta:  S=M(\mathbf{a}) \ \ \leftrightarrow
T=M(\mathbf{b})
\]
with the property that $x \in S$ avoids (123) if and only
$\Theta(x) \in T$ avoids (123).
Since any permutation of  $\mathbf{a}=(a_1,a_2, \ldots, a_k)$ can be achieved by
a sequence of transpositions of adjacent elements, it suffices to
consider the case where $\mathbf{b}$ is obtained from $\mathbf{a}$
by exchanging $a_i$ and $a_{i+1}$:
\[ \mathbf{b}=(a_1,a_2, \ldots, a_{i-1},a_{i+1},a_i,a_{i+2}, \ldots, a_k).
\]
We will do this by making use of a bijection between
permutations of $i^{a_i}(i+1)^{a_{i+1}}$ and permutations of
$i^{a_{i+1}}(i+1)^{a_{i}}$
which derives from the Greene-Kleitman symmetric
chain decomposition in the Boolean lattice \cite{GK}.

\subsection{Description of $\Theta$}

If $a_i=a_{i+1}$, the mapping $\Theta:S \rightarrow T$ is the identity.
Otherwise, if $a_i < a_{i+1}$, interchange $S$ and $T$.
If $a_i > a_{i+1}$, we proceed as follows, illustrating with the example:
\[S=(1^2)(2^1)(3^1)(4^5)(5^2)(6^7)(7^1) \rightarrow
T=(1^2)(2^1)(3^1)(4^2)(5^5)(6^7)(7^1),\]
where $T$ is  obtained from $S$ by interchanging
$a_4=5$ and $a_5=2$.

\noindent
Start with a string $x \in S$ (since $a_i > a_{i+1}$).
\begin{verbatim}
          7 5 6 6 4 6 6 4 6 6 4 6 5 3 2 4 1 1 4
\end{verbatim}
Replace $i$ by `(' and $i+1$ by `)'
\begin{verbatim}
          7 ) 6 6 ( 6 6 ( 6 6 ( 6 ) 3 2 ( 1 1 (
\end{verbatim}
Match parentheses in the usual way.
Mark unmatched left parentheses as `U'.
\begin{verbatim}
          7 ) 6 6 U 6 6 U 6 6 ( 6 ) 3 2 U 1 1 U
\end{verbatim}
Change the leftmost $a_i - a_{i+1}$ of the `U's to `)'.
Then change the remaining  `U's back to `('.
\begin{verbatim}
          7 ) 6 6 ) 6 6 ) 6 6 ( 6 ) 3 2 ) 1 1 (
\end{verbatim}
Change `)' back to $i+1$ and change `(' back to $i$.
\begin{verbatim}
          7 5 6 6 5 6 6 5 6 6 4 6 5 3 2 5 1 1 4
\end{verbatim}
This is the  string
$\Theta(x) \in T$.

\subsection{The Greene-Kleitman bijection}

A consequence of the Greene-Kleitman symmetric chain
decomposition of the Boolean lattice described in \cite{GK} is the following
bijection between the $t$-subsets and the $(n-t)$-subsets of an $n$-element
set.

If $t=n-t$ the bijection is the identity.
If $t > n-t$, define the mapping from $(n-t)$-subsets to $t$-subsets,
(i.e. from permutations of $0^{t}1^{n-t}$ to permutations of $0^{n-t}1^{t}$).
If $t < n-t$, define the  map in the reverse direction.
Assume that $t > n-t$.

For $x \in 0^{t}1^{n-t}$.
Regard `0' as `(' and `1' as `)'
and match parentheses in the usual way.
Define
\[
\tau: 0^{t}1^{n-t} \rightarrow 0^{t-1}1^{n-t+1} \]
by: $\tau(x)=y$ where $y$ is obtained from $x$ by changing
the leftmost unmatched `0' in $x$ to `1'.
So, $\tau$ can only be applied to a string with an unmatched `0'.
If $t > n-t$, there are at least $t-(n-t)=2t-n$ such.

\begin{quote}
(***): As observed in \cite{GK},
all unmatched `0's in a string are to the
right of any unmatched `1's.  Thus, the application of $\tau$ does
not change the matching.  That is, any `1' in $x$ that was matched to a `0'
is still matched in $y=\tau(x)$ to the same `0'.  So, the
unmatched `0's in $y$ are the unmatched `0's of $x$ with the leftmost
removed.  Also, changing the leftmost unmatched `0' to a `1' makes it
the rightmost unmatched `1'.
\end{quote}

\noindent
Now if $x$ has at least $j$ unmatched `0's, we could apply $\tau$
at least $j$ times.
Let $\tau^{j}(x)=\tau(\tau(...(\tau(x))...))$
        ($j$ applicatiions of $\tau$.)
It follows from (***) that
$\tau^{j}(x)=y$ where $y$ is obtained from $x$ by changing the leftmost
$j$ unmatched `0's in $x$ to `1's.
The following is the desired bijection:
\[
     \tau^{2t-n}: 0^{t}1^{n-t} \rightarrow  0^{n-t}1^{t}.
\]
To get the inverse mapping, start with a string $y \in 0^{n-t}1^{t}$.
Then $y$ has at least $t-(n-t)=2t-n$  unmatched `1's.
change the rightmost $2t-n$ unmatched `1's into `0's.

\subsection{Proof of bijection $\Theta$}

We can rephrase the mapping $\Theta: S \rightarrow T$ described in
Section 4.1 as follows:
\[
\Theta=\tau^{a_i-a_{i+1}},
\]
where
$\tau(x)=y$ and $y$ is obtained from $x$ by changing
the leftmost `$i$' in $x$ to `$i+1$'.
Then $\Theta$ is the Greene-Kleitman mapping of the preceding section,
so $\Theta$ is a bijection.  It remains to check
whether $\Theta$ preserves ``(123)-avoidance''.
We check this in steps:
\[
S \rightarrow \tau(S) \rightarrow \tau(\tau(S)) \rightarrow
\cdots \rightarrow \tau^{a_i-a_{i+1}}(S)=\Theta(S)=T
\]
and show that in each step ``(123)-avoidance" is preserved.
We show that if
$x$ is a permutation with an unmatched `$i$', then $x$ avoids (123)
iff $y=\tau(x)$ avoids (123).

Assume $x$ avoids (123) and that $x(t)=i$ is the leftmost unmatched $i$ in $x$.

\noindent
Note:

(a) If $j<t$ and $x_j=i$, then $x_j$ must be matched in $x$ to some
$i+i=x_u$ with $j<u<t$.

(b) If $j>t$ and $x_j=i+1$, then $x_j$ must be matched in $x$ to some
$i=x_u$ with $t<u<j$.

\noindent
Let $y=\tau(x)$. Then $y(t)=i+1$ and $y(i)=x(i)$ if $i \not= t$.
We show that $y$ avoids (123).
If $y$ does contains  a (123), it must involve $y_t$.
So, there exist $r,s$ such that one of the following holds:




{\bf Case (i):}  $r < s < t$ and $y_r  < y_s  < y_t =i+1$.
We cannot have $y_s  < i$, otherwise $x_r x_s x_t =y_r y_s  i $ is a
(123) in $x$.
    So, $y_s =i$.  Apply  (a) with $j=s$.
    Then  $x_r x_s x_u$ is a (123) in $x$.

{\bf Case (ii):}  $r < t < s$ and $y_r  < y_t  < y_s $.
    Note then $y_s  > i+1$.
    It can't be that $y_r  < i$, else $x_r x_t x_s =y_r  i y_s $ is a
(123) in $x$.
    So, $y_r =i$.  But then $x_r =i$.
    Apply  (a) with $j=r$.
    Then  $x_r x_u x_s$ is a (123) in $x$.

{\bf Case (iii):} $t < r < s$ and $i+1 = y_t < y_r  < y_s $.
    Then $x_t x_r x_s =i y_r y_s $ is a (123) in $x$.

Finally, we show the converse:
     if $x$ {\em does} have a (123) then so does $\tau(x)$.
Assume $x$ has a (123) and that $x_t =i$ is the leftmost unmatched $i$ in $x$.
Again $y_t =i+1$, $y_i =x_i $ if $i \not= t$.
If the (123) pattern in $x$ does not involve $x_t $, then $y=\tau(x)$ has the
same (123) pattern.  So suppose $x$ has a (123) pattern which does involve
$x_t $.
Then there exist $r,s$ such that one of the following holds:

%
%

{\bf Case (i):}  $r < s < t$ and $x_r  < x_s  < x_t =i$.
        Then $y_r y_s y_t =x_r x_s  (i+1) $ is a (123) pattern in $y$.

{\bf Case (ii):}  $r < t < s$ and $x_r  < x_t  < x_s $.
       Note then $x_r  < i$.
 If $x_s  > i+1$, then  $y_r y_t y_s =x_r  (i+1) x_s $ is a
(123) in $y$.
        Otherwise,  $x_s=i+1$.
    Apply  (b) with $j=s$.
    Then  $y_r y_u y_s$ is a (123) in $y$.

{\bf Case (iii):} $t < r < s$ and $i = x_t  < x_r  < x_s $.
    If  $x_r  > i+1$ then  $y_t y_r y_s = (i+1) x_r x_s $ is a
   (123) in $y$.
        Otherwise, $x_r =i+1$ and $x_s >i+1$.
    Apply  (b) with $j=r$.
    Then  $y_u y_r y_s$ is a (123) in $y$.
This completes the proof.

\subsection{Other patterns}

Although we can adapt $\Theta$ to show the symmetry of
$f(\mathbf{a},(321))$ in the variables $a_1, \ldots a_k$,
we note that this approach  does not work to show
symmetry for $\pi = (132)$.
For example,  the permutations in $(1^2)(2^1)(3^1)$ containing (132)
are:
\[
\{
1132,
1312,
1321\},
\]
and the permutations in $(1^1)(2^2)(3^1)$ containing (132) are:
\[
\{
2132,
1232,
1322\}.
\]
It is not possible to  get from the strings in the first set to those in
 the second simply by changing a `1' to a `2'.
The permutations (213), (231), and (312) have similar problems.
However, the same mapping $\Theta$ provides a simple bijective proof
that
$f(\mathbf{a},(12 \ldots r))$ is  symmetric in the variables $a_1, \ldots a_k$,
a result that was shown in \cite{aaa} using Schur functions and the
Robinson-Schensted-Knuth correspondence.

However the map $\Theta=\Theta(\mathbf{a},\mathbf{a'})$ can be composed with the bijections of Myers \cite{anm}, which generalize earlier constructions of Simion and Schmidt \cite{ss} to give bijective proofs of symmetry for all six patterns of three letters. Indeed, $S(\mathbf{a},\pi)$ be the set of all permutations of the multiset $1^{a_1}2^{a_2}\dots$ that avoid the pattern $\pi\in S_3$, and let
\[ SSM(\mathbf{a},\pi,\pi'): S(\mathbf{a},\pi)\rightarrow S(\mathbf{a},\pi')\]
 be the map
 of Simion-Schmidt-Myers. Then the map
\[         SSM(\mathbf{a'},(123),\pi)\circ\Theta(\mathbf{a},\mathbf{a'})\circ  SSM(\mathbf{a},\pi,(123))       \]
is a bijection between permutations of a multiset $\mathbf{a}$ that avoid $\pi$ to permutations of the multiset $\mathbf{a'}$ that avoid the same pattern $\pi$.

\section{Some open questions}
\begin{enumerate}
\item What is the asymptotic behavior of $c(n,n)$ (see
(\ref{eq:casm}) above)?
\item Investigate the equivalence classes of permutation patterns of length four
under avoidance by permutations of multisets.
\item Investigate avoidance of patterns that are not themselves permutations.
\end{enumerate}

\end{document}